# STOCHASTIC DERIVATIVES FOR FRACTIONAL DIFFUSIONS

By Sébastien Darses and Ivan Nourdin

*Université de Franche-Comté and Université Paris 6*

In this paper, we introduce some fundamental notions related to the so-called *stochastic derivatives* with respect to a given $\sigma$-field $\mathcal{Q}$. In our framework, we recall well-known results about Markov–Wiener diffusions. We then focus mainly on the case where $X$ is a fractional diffusion and where $\mathcal{Q}$ is the past, the future or the present of $X$. We treat some crucial examples and our main result is the existence of stochastic derivatives with respect to the present of $X$ when $X$ solves a stochastic differential equation driven by a fractional Brownian motion with Hurst index $H > 1/2$. We give explicit formulas.

**1. Introduction.** There exist various ways to generalize the notion of differentiation on deterministic functions. We may think of fractional derivatives or differentiation in the sense of the theory of distributions. In both cases, we lose a dynamical or geometric interpretation of tangent vectors (velocities, e.g.). In the present work, we seek to construct derivatives on stochastic processes which conserve a dynamical meaning. Our goal is motivated by the stochastic embedding of dynamical systems introduced in [2]. This procedure aims at comprehending the following question: how can we write an equation which contains the dynamical meaning of an initial ordinary differential equation and which extends this dynamical meaning to stochastic processes? We refer to [3] for more details.

Unfortunately, for most of the stochastic processes used in physical models, the limit

$$\frac{Z_{t+h} - Z_t}{h}$$

does not exist pathwise. What can be done to give a meaning to this limit? One of the main tools available is the "quantity of information" which we can use to calculate it, namely a given $\sigma$-field $\mathcal{Q}$. The idea is that one can remove the divergences which appear pathwise by averaging over a bundle of trajectories in the previous computation. This fact can be achieved by studying the behavior when $h$ goes to zero of the conditional expectation

$$\mathrm{E}\bigg[\frac{Z_{t+h} - Z_t}{h}\Big|\mathcal{Q}\bigg].$$

Such objects were introduced by Nelson in his dynamical theory of Brownian diffusion [9]. For a fixed time $t$, he calculates a forward (resp., backward) derivative with respect to a given $\sigma$-field $\mathcal{P}_t$ which can be seen as the past of the process up to time $t$ (resp., $\mathcal{F}_t$, the future of the process after time $t$). The main class with which this can be done turns out to be that of Wiener diffusions.

The purpose of this paper is, on one hand, to introduce notions which can be used to study the aforementioned quantities for general processes and, on the other hand, to treat some examples. We mainly study these notions on solutions of stochastic differential equations driven by a fractional Brownian motion with Hurst index $H \geq \frac{1}{2}$. In particular, we recall results on Wiener diffusions (case $H = \frac{1}{2}$) in our framework. We prove that for a suitable $\sigma$-algebra, the stochastic derivatives of a solution of the fractional stochastic differential equation exist and we are able to give explicit formulas.

Our paper is organized as follows. In Section 2, we recall some now classical facts on stochastic analysis for fractional Brownian motion. In Section 3, we introduce the fundamental notions related to the so-called *stochastic derivatives*. In Section 4, we study stochastic derivatives of Nelson's type for fractional diffusions. We show in Section 5 that stochastic derivatives with respect to the present turn out to be adequate tools for treating fractional Brownian motion with $H > \frac{1}{2}$. We also treat the more difficult case of a fractional diffusion.

**2. Basic notions for fractional Brownian motion.** We briefly recall some basic facts concerning stochastic calculus with respect to a fractional Brownian motion; refer to [12] for further details. Let $B = (B_t)_{t \in [0,T]}$ be a fractional Brownian motion with Hurst parameter $H \in (0,1)$ defined on a probability space $(\Omega, \mathcal{F}, \mathbb{P})$. This means that $B$ is a centered Gaussian process with the covariance function $\mathrm{E}(B_s B_t) = R_H(s,t)$, where

(1) $$R_H(s,t) = \tfrac{1}{2}(t^{2H} + s^{2H} - |t-s|^{2H}).$$

If $H = 1/2$, then $B$ is a Brownian motion. From (1), one can easily see that $\mathrm{E}|B_t - B_s|^2 = |t-s|^{2H}$, so $B$ has $\alpha$-Hölder continuous paths for any $\alpha \in (0, H)$.

2.1. *Space of deterministic integrands.* We denote by $\mathcal{E}$ the set of step $\mathbb{R}$-valued functions on $[0,T]$. Let $\mathcal{H}$ be the Hilbert space defined as the closure of $\mathcal{E}$ with respect to the scalar product

$$\langle \mathbf{1}_{[0,t]}, \mathbf{1}_{[0,s]}\rangle_{\mathcal{H}} = R_H(t,s).$$

We denote by $|\cdot|_{\mathcal{H}}$ the associated norm. The mapping $\mathbf{1}_{[0,t]} \mapsto B_t$ can be extended to an isometry between $\mathcal{H}$ and the Gaussian space $H_1(B)$ associated with $B$. We denote this isometry by $\varphi \mapsto B(\varphi)$.

When $H \in (\frac{1}{2}, 1)$, it follows from [14] that the elements of $\mathcal{H}$ may not be functions but distributions of negative order. It will be more convenient to work with a subspace of $\mathcal{H}$ which contains only functions. Such a space is the set $|\mathcal{H}|$ of all measurable functions $f$ on $[0,T]$ such that

$$|f|^2_{|\mathcal{H}|} := H(2H-1) \int_0^T \int_0^T |f(u)||f(v)||u-v|^{2H-2}\,du\,dv < \infty.$$

We know that $(|\mathcal{H}|, |\cdot|_{|\mathcal{H}|})$ is a Banach space, but that $(|\mathcal{H}|, \langle \cdot, \cdot\rangle_{\mathcal{H}})$ is not complete (see, e.g., [14]).

Moreover, we have the inclusions

(2) $$L^2([0,T]) \subset L^{1/H}([0,T]) \subset |\mathcal{H}| \subset \mathcal{H}.$$

2.2. *Fractional operators.* The covariance kernel $R_H(t,s)$ introduced in (1) can be written as

$$R_H(t,s) = \int_0^{s \wedge t} K_H(s,u)K_H(t,u)\,du,$$

where $K_H(t,s)$ is the square integrable kernel defined by

(3) $$K_H(t,s) = c_H s^{1/2-H} \int_s^t (u-s)^{H-3/2} u^{H-1/2}\,du, \qquad 0 < s < t,$$

where $c_H{}^2 = H(2H-1)\beta(2-2H, H-1/2)^{-1}$ and $\beta$ denotes the Beta function. By convention, we set $K_H(t,s) = 0$ if $s \geq t$.

Let $\mathcal{K}_H^* : \mathcal{E} \to L^2([0,T])$ be the linear operator defined by

$$\mathcal{K}_H^*(\mathbf{1}_{[0,t]}) = K_H(t,\cdot).$$

The following equality holds for any $\phi, \psi \in \mathcal{E}$:

$$\langle \phi, \psi\rangle_{\mathcal{H}} = \langle \mathcal{K}_H^*\phi, \mathcal{K}_H^*\psi\rangle_{L^2([0,T])} = \mathrm{E}(B(\phi)B(\psi)).$$

$\mathcal{K}_H^*$ then provides an isometry between the Hilbert space $\mathcal{H}$ and a closed subspace of $L^2([0,T])$.

The process $W = (W_t)_{t \in [0,T]}$ defined by

$$W_t = B((\mathcal{K}_H^*)^{-1}(\mathbf{1}_{[0,t]}))$$

is a Wiener process and the process $B$ has an integral representation of the form

$$B_t = \int_0^t K_H(t,s)\, dW_s.$$

Hence, for any $\phi \in \mathcal{H}$,

$$B(\phi) = W(\mathcal{K}_H^* \phi).$$

Let $a,b \in \mathbb{R}$, $a < b$. For any $p \geq 1$, we denote by $L^p = L^p([a,b])$ the usual Lebesgue space of functions on $[a,b]$.

Let $f \in L^1$ and $a > 0$. The left-sided and right-sided fractional Riemann–Liouville integrals of $f$ of order $\alpha$ are defined for almost all $x \in (a,b)$ by

$$I_{a+}^\alpha f(x) = \frac{1}{\Gamma(\alpha)} \int_a^x (x-y)^{\alpha-1} f(y)\, dy$$

and

$$I_{b-}^\alpha f(x) = \frac{(-1)^{-\alpha}}{\Gamma(\alpha)} \int_x^b (y-x)^{\alpha-1} f(y)\, dy,$$

respectively, where $\Gamma$ denotes the usual Euler function. These integrals extend the classical integral of $f$ when $\alpha = 1$.

If $f \in I_{a+}^\alpha(L^p)$ [resp., $f \in I_{b-}^\alpha(L^p)$] and $\alpha \in (0,1)$, then for almost all $x \in (a,b)$, the left-sided and right-sided Riemann–Liouville derivative of $f$ of order $\alpha$ are defined by

$$D_{a+}^\alpha f(x) = \frac{1}{\Gamma(1-\alpha)} \left( \frac{f(x)}{(x-a)^\alpha} + \alpha \int_a^x \frac{f(x)-f(y)}{(x-y)^{\alpha+1}}\, dy \right)$$

and

$$D_{b-}^\alpha f(x) = \frac{1}{\Gamma(1-\alpha)} \left( \frac{f(x)}{(b-x)^\alpha} + \alpha \int_x^b \frac{f(x)-f(y)}{(y-x)^{\alpha+1}}\, dy \right),$$

respectively, where $a \leq x \leq b$.

We define the operator $\mathcal{K}_H$ on $L^2([0,T])$ by

$$(\mathcal{K}_H h)(t) = \int_0^t K_H(t,s) h(s)\, ds.$$

It is an isomorphism from $L^2([0,T])$ onto $I_{0+}^{H+1/2}(L^2([0,T]))$ and it can be expressed as follows when $H > \frac{1}{2}$:

$$\mathcal{K}_H h = I_{0+}^1 s^{H-1/2} I_{0+}^{H-1/2} s^{1/2-H} h,$$

where $h \in L^2([0,T])$. The crucial point is that the functions of the space $I_{0+}^{H+1/2}(L^2([0,T]))$ are absolutely continuous when $H > \frac{1}{2}$. For these functions $\phi$, the inverse operator $\mathcal{K}_H^{-1}$ is given by

$$\mathcal{K}_H^{-1} \phi = s^{H-1/2} D_{0+}^{H-1/2} s^{1/2-H} \phi'.$$

When $H > \frac{1}{2}$, we introduce the operator $\mathcal{O}_H$ on $L^2([0,T])$ defined by

$$(4) \qquad (\mathcal{O}_H \varphi)(s) := \left(\frac{d}{dt}\mathcal{K}_H\right)(\varphi)(s) = s^{H-1/2} I_{0+}^{H-1/2} s^{1/2-H} \varphi(s).$$

Let $f:[a,b] \to \mathbb{R}$ be $\alpha$-Hölder continuous and $g:[a,b] \to \mathbb{R}$ be $\beta$-Hölder continuous with $\alpha + \beta > 1$. Then for any $s,t \in [a,b]$, the Young integral [18] $\int_s^t f \, dg$ exists and we can express it in terms of fractional derivatives (see [19]): for any $\gamma \in (1-\beta, \alpha)$, we have

$$(5) \qquad \int_s^t f \, dg = (-1)^\gamma \int_s^t D_{s+}^\gamma f(x) D_{t-}^{1-\gamma} g_{t-}(x) \, dx,$$

where $g_{t-}(x) = g(x) - g(t)$. In particular, we deduce that

$$(6) \quad \forall s < t \in [a,b] \qquad \left|\int_s^t (f(r) - f(s)) \, dg(r)\right| \leq \kappa |f|_\alpha |g|_\beta |t-s|^{\alpha+\beta},$$

where $\kappa$ is a constant depending only on $a, b, \alpha$ and $\beta$, and if $h:[a,b] \to \mathbb{R}$ and $\mu \in (0,1]$, then

$$|h|_\mu = \sup_{a \leq s < t \leq b} \frac{|h(t) - h(s)|}{|t-s|^\mu}.$$

2.3. *Malliavin calculus.* Let $\mathcal{S}$ be the set of all smooth cylindrical random variables, that is, which can be expressed as $F = f(B(\phi_1), \ldots, B(\phi_n))$ where $n \geq 1$, $f:\mathbb{R}^n \to \mathbb{R}$ is a smooth function with compact support and $\phi_i \in \mathcal{H}$. The Malliavin derivative of $F$ with respect to $B$ is the element of $L^2(\Omega, \mathcal{H})$ defined by

$$D_s^B F = \sum_{i=1}^n \frac{\partial f}{\partial x_i}(B(\phi_1), \ldots, B(\phi_n)) \phi_i(s), \qquad s \in [0,T].$$

In particular, $D_s^B B_t = \mathbf{1}_{[0,t]}(s)$. As usual, $\mathbb{D}^{1,2}$ denotes the closure of the set of smooth random variables with respect to the norm

$$\|F\|_{1,2}^2 = \mathrm{E}[F^2] + \mathrm{E}[|D_\cdot^B F|_\mathcal{H}^2].$$

The Malliavin derivative $D^B$ verifies the chain rule: if $\varphi:\mathbb{R}^n \to \mathbb{R}$ is $C_b^1$ and $(F_i)_{i=1,\ldots,n}$ is a sequence of elements of $\mathbb{D}^{1,2}$, then $\varphi(F_1, \ldots, F_n) \in \mathbb{D}^{1,2}$ and we have for any $s \in [0,T]$,

$$D_s^B \varphi(F_1, \ldots, F_n) = \sum_{i=1}^n \frac{\partial \varphi}{\partial x_i}(F_1, \ldots, F_n) D_s^B F_i.$$

The divergence operator $\delta^B$ is the adjoint of the derivative operator $D^B$. If a random variable $u \in L^2(\Omega, \mathcal{H})$ belongs to the domain of the divergence operator, that is, if it verifies

$$|\mathrm{E}\langle D^B F, u\rangle_\mathcal{H}| \leq c_u \|F\|_{L^2} \qquad \text{for any } F \in \mathcal{S},$$

then $\delta^B(u)$ is defined by the duality relationship

$$E(F\delta^B(u)) = E\langle D^B F, u\rangle_{\mathcal{H}}$$

for every $F \in \mathbb{D}^{1,2}$.

2.4. *Pathwise integration with respect to B.* If $X = (X_t)_{t\in[0,T]}$ and $Z = (Z_t)_{t\in[0,T]}$ are two continuous processes, then we define the forward integral of $Z$ with respect to $X$, in the sense of Russo–Vallois, by

$$(7) \qquad \int_0^{\cdot} Z_s \, dX_s = \lim_{\varepsilon \to 0} \text{ucp}\, \varepsilon^{-1} \int_0^{\cdot} Z_s(X_{s+\varepsilon} - X_s)\, ds, \qquad t \in [0, T],$$

provided the limit exists. Here "ucp" means "uniform convergence in probability." If $X$ (resp., $Z$) has a.s. Hölder continuous paths of order $\alpha$ (resp., $\beta$) with $\alpha + \beta > 1$, then $\int_0^{\cdot} Z_s \, dX_s$ exists and coincides with the usual Young integral: see [15], Proposition 2.12.

2.5. *Stochastic differential equation driven by B.* Here, we assume that $H > 1/2$. We denote by $C_b^k$ the set of all functions whose derivatives from order 1 to order $k$ are bounded. If $\sigma \in C_b^2$ and $b \in C_b^1$, then the equation

$$(8) \qquad X_t = x_0 + \int_0^t \sigma(X_s)\, dB_s + \int_0^t b(X_s)\, ds, \qquad t \in [0, T],$$

admits a unique solution $X$ in the set of processes whose paths are Hölder continuous of order $\alpha > 1 - H$. Here, the integral with respect to $B$ is in the sense of Russo–Vallois; see (7). Moreover, we have a Doss–Sussmann-type [5, 16] representation,

$$X_t = \phi(A_t, B_t), \qquad t \in [0, T],$$

where $\phi$ and $A$ are given, respectively, by

$$\frac{\partial \phi}{\partial x_2}(x_1, x_2) = \sigma(\phi(x_1, x_2)), \qquad \phi(x_1, 0) = x_1, \ x_1, x_2 \in \mathbb{R}$$

and

$$A'_t = \exp\left(-\int_0^{B_t} \sigma'(\phi(A_t, s))\, ds\right) b(\phi(A_t, B_t)), \qquad A_0 = x_0, \ t \in [0, T].$$

Using this representation, we can show that $X$ belongs to $\mathbb{D}^{1,2}$ and that

$$D_s^B X_t = \sigma(X_s) \exp\left(\int_s^t b'(X_u)\, du + \int_s^t \sigma'(X_u)\, dB_u\right) \mathbf{1}_{[0,t]}(s), \qquad s, t \in [0, T];$$

see [10], proof of Theorem B.

**3. Notions related to stochastic derivatives.** Let $(Z_t)_{t \in [0,T]}$ be a stochastic process defined on $(\Omega, \mathcal{F}, \mathbb{P})$. In the sequel, we always assume that for any $t \in [0,T]$, $Z_t \in L^2(\Omega, \mathcal{F}, \mathbb{P})$. For all $t \in (0,T)$ and $h \neq 0$ such that $t + h \in (0,T)$, we set

$$\Delta_h Z_t = \frac{Z_{t+h} - Z_t}{h}.$$

3.1. *Stochastic derivatives in a strong sense.*

DEFINITION 1. Set $t \in (0,T)$. We say that $\mathcal{A}^t$ (resp., $\mathcal{B}^t$) is a *forward differentiating $\sigma$-field* (resp., *backward differentiating $\sigma$-field*) for $Z$ at $t$ if $\mathrm{E}[\Delta_h Z_t | \mathcal{A}^t]$ (resp., $\mathrm{E}[\Delta_{-h} Z_t | \mathcal{B}^t]$) converges in probability when $h \downarrow 0^+$. In these cases, we define the so-called forward and backward derivatives,

(9) $$D_+^{\mathcal{A}^t} Z_t = \lim_{h \downarrow 0^+} \mathrm{E}[\Delta_h Z_t | \mathcal{A}^t],$$

(10) $$D_-^{\mathcal{B}^t} Z_t = \lim_{h \downarrow 0^+} \mathrm{E}[\Delta_{-h} Z_t | \mathcal{B}^t].$$

The set of all forward (resp., backward) differentiating $\sigma$-fields for $Z$ at time $t$ is denoted by $\mathcal{M}_Z^{+(t)}$ (resp., $\mathcal{M}_Z^{-(t)}$). The intuition we can have is that the larger $\mathcal{M}^{\pm(t)}$ is, the more regular $Z$ is at time $t$. For instance, one obviously has that $\{\varnothing, \Omega\} \in \mathcal{M}_Z^{+(t)}$ (resp., $\in \mathcal{M}_Z^{-(t)}$) if and only if $s \mapsto \mathrm{E}(Z_s)$ is right differentiable (resp., left differentiable) at time $t$. At the opposite extreme, one has that $\mathcal{F} \in \mathcal{M}_Z^{+(t)}$ (resp., $\in \mathcal{M}_Z^{-(t)}$) if and only if $s \mapsto Z_s$ is a.s. right differentiable (resp., left differentiable) at time $t$.

DEFINITION 2. We say that $(\mathcal{A}^t, \mathcal{B}^t)_{t \in (0,T)}$ is a *differentiating collection of $\sigma$-fields* for $Z$ if for any $t \in (0,T)$, $\mathcal{A}^t$ (resp., $\mathcal{B}^t$) is a forward (resp., backward) differentiating $\sigma$-field for $Z$ at $t$. If $\mathcal{A}^t = \mathcal{B}^t$ for any $t \in (0,T)$, we write, for simplicity, $(\mathcal{A}^t)_{t \in (0,T)}$ instead of $(\mathcal{A}^t, \mathcal{B}^t)_{t \in (0,T)}$.

On one hand, we may introduce the following definition.

DEFINITION 3. Set $t \in (0,T)$. We say that $\mathcal{A}^t$ (resp. $\mathcal{B}^t$) is a *nondegenerate forward $\sigma$-field* (resp., *nondegenerate backward $\sigma$-field*) for $Z$ at $t$ if it is forward (resp., backward) differentiating at $t$ and if

(11)  for any $c \in \mathbb{R}$   $\mathrm{P}(D_+^{\mathcal{A}^t} Z_t = c) < 1$   [resp., $\mathrm{P}(D_-^{\mathcal{B}^t} Z_t = c) < 1$].

For instance, if $Z$ is a process such that $s \mapsto \mathrm{E}(Z_s)$ is differentiable at $t \in (0,T)$, then $\{\varnothing, \Omega\}$ is a forward and backward differentiating $\sigma$-field at $t$ but is degenerate. Let us also note that condition (11) is obviously equivalent to $\mathrm{Var}(D_+^{\mathcal{A}^t} Z_t) \neq 0$ [resp., $\mathrm{Var}(D_-^{\mathcal{B}^t} Z_t) \neq 0$] when $D_+^{\mathcal{A}^t} Z_t \in L^2(\Omega)$ [resp., $D_-^{\mathcal{B}^t} Z_t \in L^2(\Omega)$].

On the other hand, one could hope that such stochastic derivatives conserve the property which holds for ordinary derivatives on functions: "it can discriminate the constants among the other processes." So we introduce the following.

DEFINITION 4.  We say that $(\mathcal{A}^t, \mathcal{B}^t)_{t \in (0,T)}$ is *a discriminating collection of $\sigma$-fields* for $Z$ if $(\mathcal{A}^t, \mathcal{B}^t)_{t \in (0,T)}$ is a differentiating collection of $\sigma$-fields for $Z$ and if it satisfies the following property:

$$(\forall t \in (0,T),\ D_+^{\mathcal{A}^t} Z_t = D_-^{\mathcal{B}^t} Z_t = 0) \Rightarrow Z \quad \text{is a.s. a constant process on } [0, T].$$

As in Definition 2, we write, for simplicity, $(\mathcal{A}^t)_{t \in (0,T)}$ instead of $(\mathcal{A}^t, \mathcal{B}^t)_{t \in (0,T)}$ when $\mathcal{A}^t = \mathcal{B}^t$ for any $t \in (0,T)$.

An obvious example of a discriminating collection of $\sigma$-fields for a process with differentiable paths is $\{\mathcal{A}^t = \mathcal{F},\ t \in (0,T)\}$. If $Z$ is a process such that $s \mapsto \mathrm{E}(Z_s)$ is differentiable on $(0,T)$, then the collection $\{\mathcal{A}^t = \{\varnothing, \Omega\}, t \in [0,T]\}$ is differentiating, but, in general, not discriminating.

Let us now consider a more advanced example. Let $B = (B_t)_{t \in [0,T]}$ be a fractional Brownian motion with Hurst index $H \in (1/2, 1)$. Let us denote by $\mathcal{P}_t$ the $\sigma$-field generated by $B_s$ for $0 \leq s \leq t$ and, if $g : \mathbb{R} \to \mathbb{R}$, by $\mathcal{T}_t^g$ the $\sigma$-field generated by $g(B_t)$.

PROPOSITION 5.  *1. For any $t \in (0,T)$, $\mathcal{P}_t$ is not a forward differentiating $\sigma$-field for $B$ at $t$.*

*2. For any even function $g : \mathbb{R} \to \mathbb{R}$ and for any $t \in (0,T)$, $\mathcal{T}_t^g$ is a forward and backward differentiating (but generate) $\sigma$-field for $B$ at $t$.*

*3. For any $t \in (0,T)$, $\mathcal{T}_t^{\mathrm{id}}$ is a forward and backward differentiating and nondegenerate $\sigma$-field for $B$ at $t$.*

PROOF.  1. We refer to Proposition 10 of [4] for a proof. This result is extended to the case of Volterra processes in the current paper; see Proposition 13.

2. Since $B$ and $-B$ have the same law, it follows that $\mathrm{E}[\Delta_h B_t | g(B_t)] = 0$ for any $t \in (0,T)$ and $h \neq 0$ such that $t + h \in (0,T)$. The conclusion follows easily.

3. Using a linear Gaussian regression we can write

$$\mathrm{E}[\Delta_h B_t | B_t] = \frac{(1 + h/t)^{2H} - 1 - (h/t)^{2H}}{2} B_t \xrightarrow[h \to 0]{} H \frac{B_t}{t} \quad \text{in probability.}$$

Since $\mathrm{Var}(Ht^{-1}B_t) > 0$, $\mathcal{T}_t^{\mathrm{id}}$ is nondegenerate. □

Thus, for fractional Brownian motion, stochastic derivatives with respect to the present (i.e., with respect to $\mathcal{T}_t^{\mathrm{id}}$) turn out to be adequate tools (see Section 5 below, for a more precise study).

3.2. *Stochastic derivatives in a weak sense.* A way to weaken Definition 1 is to consider stochastic derivatives as follows.

DEFINITION 6. Set $t \in (0, T)$ and let $\mathcal{A}$ be a sub-$\sigma$-field of $\mathcal{F}$. We say that $Z$ is weak forward differentiable with respect to $\mathcal{A}$ at $t$ if $\lim_{h \downarrow 0^+} \mathrm{E}[V \Delta_h Z_t]$ exists, for all random variables $V$ belonging to a dense subspace of the closed subspace $L^2(\Omega, \mathcal{A}, \mathbb{P}) \subset L^2(\Omega, \mathcal{F}, \mathbb{P})$.

We similarly define the notion of *weak backward differentiation* with respect to $\mathcal{A}$ at $t$ by considering $\Delta_{-h} Z_t$ instead of $\Delta_h Z_t$.

If the process $Z$ is weak forward differentiable at $t$ and such that the sequence $(\Delta_h Z_t)_h$ is bounded in $L^2(\Omega)$, then we can associate with it a weak forward stochastic derivative with respect to $\mathcal{A}$ at $t$. Indeed, in that case, let us denote by $\Theta$ the dense subspace involved. The linear form $\psi : V \mapsto \lim_{h \downarrow 0^+} \mathrm{E}[V \Delta_h Z_t]$ defined on $\Theta \subset L^2(\Omega, \mathcal{A}, \mathbb{P})$ is continuous and so can be extended in a unique continuous linear form on $L^2(\Omega, \mathcal{A}, \mathbb{P})$, still denoted by $\psi$. Thus, there exists a unique $Z'_t \in L^2(\Omega, \mathcal{A}, \mathbb{P})$ such that $\psi(V) = \mathrm{E}[Z'_t V]$. One can easily show that $Z'_t$ does not depend on $\Theta$. We will say that $Z'_t$ is the *weak forward stochastic derivative* of $Z$ with respect to $\mathcal{A}$ at $t$.

REMARK 7. The boundedness of $(\Delta_h Z_t)_h$ in $L^2(\Omega)$ may appear to be quite a restrictive condition (e.g., it is not satisfied for a fractional Brownian motion). But it allows us to relate our notion to the usual notion of weak limit.

If $\mathcal{A}^t$ (resp., $\mathcal{B}^t$) is a forward (resp., backward) differentiating $\sigma$-field for $Z$ at $t$ and if, moreover, the convergence (9) [resp., (10)] also holds in $L^2$, then $Z$ is weak forward (resp., backward) differentiable with respect to $\mathcal{A}^t$ (resp., $\mathcal{B}^t$) at $t$. But the converse is not true in general.

Let $\Upsilon$ be the set of the so-called fractional diffusions $X = (X_t)_{t \in [0,T]}$ defined by

$$(12) \qquad X_t = x_0 + \int_0^t \sigma_s \, dB_s + \int_0^t b_s \, ds, \qquad t \in [0, T].$$

Here, $\sigma$ and $b$ are processes which are adapted with respect to the natural filtration associated with $B$ and $X$ and satisfying the following conditions: $\sigma \in C^\alpha$ a.s. with $\alpha > 1 - H$ and $b \in L^1([0, T])$ a.s.

LEMMA 8. *The decomposition (12) is unique. That is, if*

$$x_0 + \int_0^t \sigma_s \, dB_s + \int_0^t b_s \, ds = \tilde{x}_0 + \int_0^t \tilde{\sigma}_s \, dB_s + \int_0^t \tilde{b}_s \, ds, \qquad t \in [0,T], \tag{13}$$

*then $x_0 = \tilde{x}_0$, $\sigma = \tilde{\sigma}$ and $b = \tilde{b}$.*

PROOF. The equality $x_0 = \tilde{x}_0$ is obvious and (13) is then equivalent to

$$\int_0^t (\sigma_s - \tilde{\sigma}_s) \, dB_s = \int_0^t (\tilde{b}_s - b_s) \, ds, \qquad t \in [0,T],$$

which implies, by setting $t_k = \frac{kT}{n}$, that

$$(|\sigma_{t_k} - \tilde{\sigma}_{t_k}||B_{t_{k+1}} - B_{t_k}|)^{1/H}$$

$$= \left| \int_{t_k}^{t_{k+1}} (b_s - \tilde{b}_s) \, ds + \int_{t_k}^{t_{k+1}} (\sigma_s - \sigma_{t_k}) \, dB_s + \int_{t_k}^{t_{k+1}} (\tilde{\sigma}_s - \tilde{\sigma}_{t_k}) \, dB_s \right|^{1/H}$$

$$\leq C \left[ \left| \int_{t_k}^{t_{k+1}} (b_s - \tilde{b}_s) \, ds \right|^{1/H} + \left| \int_{t_k}^{t_{k+1}} (\sigma_s - \sigma_{t_k}) \, dB_s \right|^{1/H} \right.$$

$$\left. + \left| \int_{t_k}^{t_{k+1}} (\tilde{\sigma}_s - \tilde{\sigma}_{t_k}) \, dB_s \right|^{1/H} \right].$$

We easily deduce, using (6), that

$$\lim_{n \to \infty} \sum_{k=0}^{n-1} |\sigma_{t_k} - \tilde{\sigma}_{t_k}|^{1/H} |B_{t_{k+1}} - B_{t_k}|^{1/H} = 0 \qquad \text{in probability.}$$

But, on the other hand, it is easy to obtain (see, e.g., Theorem 4.4 in [7]) that

$$\lim_{n \to \infty} \sum_{k=0}^{n-1} |\sigma_{t_k} - \tilde{\sigma}_{t_k}|^{1/H} |B_{t_{k+1}} - B_{t_k}|^{1/H} = \int_0^T |\sigma_s - \tilde{\sigma}_s|^{1/H} \, ds \qquad \text{in probability.}$$

We deduce that $\sigma = \tilde{\sigma}$ and so $b = \tilde{b}$. □

In Section 4, we will see that the past of $X \in \Upsilon$ before time $t$ is not, in general, a forward differentiating $\sigma$-field at time $t$. We will see in Section 5 that the present of $X \in \Upsilon$ is, in general, a differentiating collection of $\sigma$-fields.

However, $X$ is weak differentiable for a large class of $\sigma$-fields. We introduce the set $\mathcal{S}^b$ of all r.v.'s $\varphi(B(\phi_1), \ldots, B(\phi_n)) \in \mathcal{S}$ such that $\phi_1, \ldots, \phi_n$ are bounded functions.

Let $\wp$ be the set of all sub-$\sigma$-fields $\mathcal{A} \subset \mathcal{F}$ such that $L^2(\Omega, \mathcal{A}, \mathbb{P}) \cap \mathcal{S}^b$ is dense in $L^2(\Omega, \mathcal{A}, \mathbb{P})$. For instance, any $\sigma$-field can be expressed as $\mathcal{A}^{[r,s]} = \varsigma(B_v, r \leq v \leq s)$ belongs to $\wp$ (see, e.g., [12], page 24).

PROPOSITION 9. *Let $\mathcal{A} \in \wp$ and $t \in (0,T)$. Let $X \in \Upsilon$ be given by (12), satisfying the following conditions:*

(i) *the map $s \mapsto b_s$ is continuous from $(0,T)$ into $L^1(\Omega)$;*
(ii) *for all $s \in [0,T]$, $\sigma_s \in \mathbb{D}^{1,2}$ and $\sup_{s \in [0,T]} \mathrm{E}|D_s^B \sigma_t| < +\infty$;*
(iii) *$\mathrm{E}(|\sigma|_\alpha^p) < +\infty$ for some $p > 1$ and $\alpha > 1 - H$.*

*Then $X$ is weak forward and backward differentiable at $t$ with respect to $\mathcal{A}$.*

PROOF. For simplicity, we only prove the forward case, the backward case being similar. Let $t \in (0,T)$.

We write

$$(14) \quad X_{t+h} - X_t = \sigma_t(B_{t+h} - B_t) + \int_t^{t+h} b_s \, ds + \int_t^{t+h} (\sigma_s - \sigma_t) \, dB_s.$$

First, we treat the second term of the right-hand side of (14). Let $V \in L^2(\Omega, \mathcal{A}, \mathbb{P}) \cap \mathcal{S}^b$. Since $V$ is bounded and the map $s \mapsto b_s$ is continuous from $(0,T)$ into $L^1(\Omega)$, the function $s \mapsto \mathrm{E}[Vb_s]$ is continuous. We then deduce, by means of Fubini's theorem, that

$$(15) \quad \lim_{h \downarrow 0} \frac{1}{h} \mathrm{E}\left[V \int_t^{t+h} b_s \, ds\right] ds = \mathrm{E}[Vb_t].$$

Then using inequality (6) and the hypothesis $\mathrm{E}(|\sigma|_\alpha^p) < +\infty$, the following limit holds:

$$(16) \quad \lim_{h \to 0} \frac{1}{h} \mathrm{E}\left[V \int_t^{t+h} (\sigma_s - \sigma_t) \, dB_s\right] = 0.$$

Finally, we show that the limit

$$\lim_{h \downarrow 0} \mathrm{E}[\sigma_t V \Delta_h B_t]$$

exists. Since $\sigma_t V \in \mathbb{D}^{1,2}$ (see Exercise 1.2.13 in [11]), we have

$$\mathrm{E}[\sigma_t(B_{t+h} - B_t)V] = \mathrm{E}[\delta^B(\mathbf{1}_{[t,t+h]})\sigma_t V]$$
$$= \mathrm{E}[\sigma_t \langle \mathbf{1}_{[t,t+h]}, D^B V \rangle_{\mathcal{H}}] + \mathrm{E}[V \langle \mathbf{1}_{[t,t+h]}, D^B \sigma_t \rangle_{\mathcal{H}}].$$

Condition (ii) and the fact that $V \in \mathcal{S}^b$ allow us, in particular, to write

$$(17) \quad \mathrm{E}[\sigma_t(B_{t+h} - B_t)V] = H(2H-1)(\Psi_{t,h}(\sigma_t, V) + \Psi_{t,h}(V, \sigma_t)),$$

where

$$\Psi_{t,h}(X,Y) = \mathrm{E}\left[X \int_0^T D_s^B Y \int_t^{t+h} |v - s|^{2H-2} \, dv \, ds\right].$$

When $X$ or $Y$ denotes $\sigma_t$ or $V$, Fubini's theorem yields

$$\Psi_{t,h}(X,Y) = \int_t^{t+h} f(v, X, Y) \, dv$$

with

$$f(v, X, Y) = \int_0^T \mathrm{E}[X D_s^B Y]|v - s|^{2H-2}\, ds.$$

We have, due to condition (ii) and the fact that $V \in \mathcal{S}^b$, that

$$|f(v, X, Y) - f(w, X, Y)| \leq C(X, Y) \int_0^T ||v - s|^{2H-2} - |w - s|^{2H-2}|\, ds,$$

where $C(X, Y)$ is a constant depending only on $X$ and $Y$.

The previous integral tends to 0 when $w$ tends to $v$.
The continuity of the function $v \mapsto f(v, X, Y)$ follows.
Therefore, the limit

$$\lim_{h \to 0} h^{-1} \mathrm{E}[\sigma_t(B_{t+h} - B_t)V]$$

exists and equals

$$H(2H - 1)\mathrm{E}\left[\sigma_t \int_0^T D_s^B V |t - s|^{2H-2}\, ds + V \int_0^T D_s^B \sigma_t |t - s|^{2H-2}\, ds\right]. \quad \square$$

**4. Stochastic derivatives of Nelson's type.** Let $Z$ be a stochastic process defined on $(\Omega, \mathcal{F}, \mathbb{P})$. We define the past of $Z$ before time $t$,

$$\mathcal{P}_t^Z := \varsigma(Z_s, 0 \leq s \leq t)$$

and the future of $Z$ after time $t$,

$$\mathcal{F}_t^Z := \varsigma(Z_s, t \leq s \leq T).$$

If $\mathcal{P}_t^Z$ and $\mathcal{F}_t^Z$ are, respectively, forward and backward differentiating $\sigma$-fields for $Z$ at $t$, we call $D_+^{\mathcal{P}_t^Z} Z_t$ and $D_-^{\mathcal{F}_t^Z} Z_t$, respectively, the forward and backward stochastic derivatives *of Nelson's type* in reference to Nelson's work [9]. In the sequel, we denote them by $D_+^{\mathcal{P}} Z_t$ and $D_-^{\mathcal{F}} Z_t$, for simplicity.

4.1. *The case of Wiener diffusions.* We denote by $\Lambda$ the space of diffusion processes $X$ satisfying the following conditions.

1. $X$ solves the stochastic differential equation

(18) $$dX_t = b(t, X_t)\, dt + \sigma(t, X_t)\, dW_t, \qquad X_0 = x_0,$$

where $x_0 \in \mathbb{R}^d$, $b : [0, T] \times \mathbb{R}^d \to \mathbb{R}^d$ and $\sigma : [0, T] \times \mathbb{R}^d \to \mathbb{R}^d \otimes \mathbb{R}^d$ are Borel measurable functions satisfying the following hypothesis: there exists a constant $K > 0$ such that for every $x, y \in \mathbb{R}^d$, we have

$$\sup_t(|\sigma(t, x) - \sigma(t, y)| + |b(t, x) - b(t, y)|) \leq K|x - y|,$$

$$\sup_t(|\sigma(t, x)| + |b(t, x)|) \leq K(1 + |x|).$$

2. For any $t \in (0,T)$, $X_t$ has a density $p_t$.
3. Setting $a_{ij} = (\sigma\sigma^*)_{ij}$, for any $i,j \in \{1,\ldots,n\}$ any $t_0 \in (0,T)$ and any bounded open set $O \subset \mathbb{R}^d$,

$$\int_{t_0}^{T} \int_{O} |\partial_j(a_{ij}(t,x)p_t(x))| \, dx \, dt < +\infty.$$

4. The functions $b$ and $(t,x) \mapsto \frac{1}{p_t(x)} \partial_j(a_{ij}(t,x)p_t(x))$ are bounded, belong to $C^{1,2}([0,T] \times \mathbb{R}^d)$ and have all first- and second-order derivatives bounded [we use the usual convention that the term involving $\frac{1}{p_t(x)}$ is 0 if $p_t(x) = 0$].

These conditions are introduced in [8] and ensure the existence of a drift and a diffusion coefficient for the time-reversed process $\overline{X}_t := X_{T-t}$. Föllmer focuses in [6], Proposition 2.5, on the important relation between drifts and derivatives of Nelson's type. It allows the computation the drift of the time reversal of a Brownian diffusion with constant diffusion coefficient, both in the Markov and non-Markov case (see Theorem 3.10 and 4.7 in [6]).

For a Markov diffusion with a rather general diffusion coefficient, we have the following result.

THEOREM 10. *Let $X \in \Lambda$ be given by (18). Then $X$ is a Markov diffusion with respect to $\mathcal{P}^X$ and $\mathcal{F}^X$. Moreover, $\mathcal{P}^X$ and $\mathcal{F}^X$ are differentiating and, in general, nondegenerate:*

$$D_+^{\mathcal{P}} X_t = b(t, X_t),$$

$$(D_-^{\mathcal{F}} X_t)_i = b_i(t, X_t) - \frac{1}{p_t(X_t)} \sum_j \partial_j(a_{ij}(t,X_t)p_t(X_t)),$$

*where we use the convention that the term involving $\frac{1}{p_t(x)}$ is 0 if $p_t(x) = 0$.*

We refer to [3] for a proof; it is based on the proof of Proposition 4.1 in [17] and Theorem 2.3 in [8].

4.2. *The case of fractional Brownian motion and Volterra processes.* Let $K$ be an $L^2$-kernel, that is, a function $K:[0,T] \times [0,T] \to \mathbb{R}$ satisfying $\int_{[0,T]^2} K(t,s)^2 \, dt \, ds < +\infty$. We denote by $\frac{\partial^+ K}{\partial t}$ the *right* derivative of $K$ with respect to $t$ (with the convention that it equals $+\infty$ if it does not exist).

We assume, moreover, that $K$ is Volterra, that is, that it vanishes on $\{(t,s) \in [0,T]^2 : s > t\}$ and is nondegenerate. In other words, the family $\{K(t,\cdot), t \in [0,T]\}$ is free and spans a vector space dense in $L^2([0,T])$. With such a kernel $K$, we associate the so-called *Volterra process*

(19) $$G_t = \int_0^t K(t,s) \, dW_s, \qquad 0 \le t \le T,$$

where $W$ denotes a standard Brownian motion. The assumptions made on $K$ imply, in particular, that the natural filtrations associated with $W$ and $G$ are the same (see, e.g., [1], Remark 3).

PROPOSITION 11. *Let $t \in (0,T)$ and $G$ be a Volterra process associated with a nondegenerate Volterra kernel $K$ satisfying the condition*

$$
(20) \qquad \frac{K(t+h,\cdot) - K(t,\cdot)}{h} \xrightarrow[h\downarrow 0]{} \frac{\partial^+ K}{\partial t} \qquad \text{in } L^2([0,t]).
$$

*The forward Nelson derivative $D^{\mathcal{P}}_+ G_t$ at $t$ exists if and only if $\int_0^t \frac{\partial^+ K}{\partial t}(t,s)^2\,ds < +\infty$. In this case, we have $D^{\mathcal{P}}_+ G_t = \int_0^t \frac{\partial^+ K}{\partial t}(t,s)\,dW_s$ and $\mathcal{P}^G_t$ is nondegenerate at $t$ if and only if $\int_0^t \frac{\partial^+ K}{\partial t}(t,s)^2\,ds > 0$.*

PROOF. We adapt the proof of [4], Proposition 10. Using the representation (19), we deduce that

$$
\begin{aligned}
\mathrm{E}[\Delta_h G_t | \mathcal{P}^G_t] &= \mathrm{E}[\Delta_h G_t | \mathcal{P}^W_t] \\
&= \frac{1}{h} \int_0^t [K(t+h,s) - K(t,s)]\,dW_s =: Z_h.
\end{aligned}
$$

Note that $Z = (Z_h)_{h>0}$ is a centered Gaussian process. First, assume that $\int_0^t \frac{\partial^+ K}{\partial t}(t,s)^2\,ds = +\infty$. It is a classical result that if $Z_h$ converges in probability as $h \downarrow 0$, then $\mathrm{Var}(Z_h)$ converges as $h \downarrow 0$. But, from Fatou's lemma, we deduce that

$$
\liminf_{h\downarrow 0} \mathrm{Var}(Z_h) \geq \int_0^t \frac{\partial^+ K}{\partial t}(t,s)^2\,ds = +\infty.
$$

Thus, $Z_h$ does not converge in probability as $h \downarrow 0$. Conversely, assume that $\int_0^t \frac{\partial^+ K}{\partial t}(t,s)^2\,ds < +\infty$. In this case, assumption (20) implies that $Z_h \to \int_0^t \frac{\partial^+ K}{\partial t}(t,s)\,dW_s$ in probability as $h \downarrow 0$. In other words, $D^{\mathcal{P}}_+ G_t$ exists and equals $\int_0^t \frac{\partial^+ K}{\partial t}(t,s)\,dW_s$. We easily deduce that $\mathcal{P}^G_t$ is nondegenerate at $t$ if and only if $\int_0^t \frac{\partial^+ K}{\partial t}(t,s)^2\,ds > 0$. □

The result of Proposition 10 in [4] is then a particular case: if $B$ denotes a fractional Brownian motion with Hurst index $H \in (0,1/2) \cup (1/2,1)$ and if $t \in (0,T)$, then $D^{\mathcal{P}}_+ B_t$ does not exist. Indeed, we have $B_t = \int_0^t K_H(t,s)\,dW_s$, where $K_H$ is the nondegenerate Volterra kernel given by (3) and verifying

$$
\frac{\partial K_H}{\partial t}(t,s) = c_H \left(\frac{t}{s}\right)^{H-1/2}(t-s)^{H-3/2}.
$$

REMARK 12. For a stochastic process $Z$, let us define

(21) $$\xi(Z) = \text{Leb}\{t \in [0,T],\ D_+^{\mathcal{P}} Z_t \text{ exists}\}.$$

For instance, if $B$ is a fractional Brownian motion with Hurst index $H \in (0,1)$, then $\xi(B) = T$ if $H = 1/2$ and $\xi(B) = 0$ otherwise. A real $c \in [0,T]$ being fixed, it is, in fact, not difficult, using Proposition 11, to construct a continuous process $Z$ such that $\xi(Z) = c$. For instance, we can consider the Volterra process associated with the Volterra kernel

$$K(t,s) = \begin{cases} (t-s)^{H(t)}, & \text{if } s \le t, \\ 0, & \text{otherwise,} \end{cases}$$

$$\text{with } H(t) = \begin{cases} 0, & \text{if } t \le c, \\ (t-c) \wedge 1/4, & \text{if } t > c. \end{cases}$$

The study of backward derivatives seems to be more difficult. Among the difficulties, we note that it is not easy to obtain backward representation of fractional diffusions (see [4]). However, for a fractional Brownian motion, we are able to prove the following proposition.

PROPOSITION 13. *Set $H > 1/2$. The limit*

$$\lim_{h \downarrow 0} \text{E}\left[ \frac{B_t - B_{t-h}}{h} \Big| \mathcal{F}_t^B \right]$$

*exists neither as an element of $L^p(\Omega)$ for any $p \in [1,\infty)$ nor as an almost sure limit.*

PROOF. Fix $t \in (0,T)$ and set

$$G_h := \text{E}\left[ \frac{B_t - B_{t-h}}{h} \Big| \mathcal{F}_t^B \right] \quad \text{and} \quad Z_h := \text{E}\left[ \frac{B_t - B_{t-h}}{h} \Big| B_t, B_{t+h} \right].$$

Since $(G_h)_{h>0}$ is a family of Gaussian random variables, it suffices to prove that $\text{Var}(G_h)$ diverges when $h$ goes to 0.

We have $Z_h = \text{E}[G_h | B_t, B_{t+h}]$. So, by Jensen's inequality, $Z_h^2 \le \text{E}[G_h^2 | B_t, B_{t+h}]$ and $\text{Var}(Z_h) \le \text{Var}(G_h)$. Let us show that $\lim_{h \downarrow 0} \text{Var}(Z_h) = +\infty$.

The covariance matrix of the Gaussian vector $(B_{t-h} - B_t, B_t, B_{t+h})$ is

$$\begin{pmatrix} a & v \\ v^* & M \end{pmatrix},$$

where $a = \text{Var}(B_{t-h} - B_t)$, $v = (R(t-h, t) - R(t,t); R(t-h, t+h) - R(t, t+h))$ and

$$M = \begin{pmatrix} R(t,t) & R(t, t+h) \\ R(t, t+h) & R(t+h, t+h) \end{pmatrix}.$$

Since $d_h := R(t,t)R(t+h,t+h) - R(t+h,t)^2 \neq 0$, $M$ is invertible. Therefore, $hZ_h = vM^{-1}Q^*$, where $Q = (B_t, B_{t+h})$. Since $M = \mathrm{E}[Q^*Q]$, we deduce that

$$\mathrm{Var}(hZ_h) = \mathrm{E}[vM^{-1}Q^*(vM^{-1}Q^*)^*] = vM^{-1}v^*.$$

Hence,

$$\mathrm{Var}(hZ_h) = \frac{1}{d_h}(R(t+h,t+h)v_1^2 - 2R(t+h,t)v_1v_2 + R(t,t)v_2^2).$$

This expression is homogeneous in $t^{2H}$, so we henceforth work with $t=1$. Tedious computation gives $d_h \sim h^{2H}$ as $h \downarrow 0$. Moreover, we note that $v_2 = v_1 + \mathrm{c}h^{2H}$, where c is a constant depending only on $H$. Thus,

$$d_h \mathrm{Var}(hZ_h) = v_1 \mathrm{c} h^{2H}(1 - (1+h)^{2H} + h^{2H}) + h^{2H}v_1^2 + \mathrm{c}^2 h^{4h}.$$

Since $2H > 1$ and the function $x \mapsto x^{2H}$ is derivable, the quantities $\frac{v_1}{h}$ and $\frac{1-(1+h)^{2H}+h^{2H}}{h}$ converge as $h \downarrow 0$. But $2H < 2$ and $\frac{h^{4H}}{h^2 h^{2H}} = h^{2H-2} \to +\infty$ as $h \downarrow 0$. Thus,

$$\lim_{h \downarrow 0} \mathrm{Var}(Z_h) = +\infty,$$

which concludes the proof. $\square$

4.3. *The case of fractional diffusions.*

PROPOSITION 14. *Let $X \in \Upsilon$ be given by (12) and satisfy the following conditions:* $\mathrm{E}(\int_0^T |b_s|\,ds) < +\infty$ *and* $\mathrm{E}(|\sigma|_\alpha^p) < +\infty$ *for some $p > 1$ and $\alpha > 1 - H$. If for any $t \in (0,T)$, $\sigma_t \neq 0$ a.s., then for almost all $t \in (0,T)$, $\mathcal{P}_t^X$ is not a forward differentiating $\sigma$-field for $X$ at $t$.*

PROOF. Remember we assumed that $\sigma$ and $b$ are adapted with respect to the natural filtration associated with $B$ and $X$, see (12). In particular, we deduce from (12) that $\mathcal{P}_t^X \subset \mathcal{P}_t^B$. Since we can also write

$$B_t = \int_0^t \frac{1}{\sigma_s}\,dX_s - \int_0^t \frac{b_s}{\sigma_s}\,ds,$$

we finally have $\mathcal{P}_t^X = \mathcal{P}_t^B$.

Thus, we deduce that $\mathrm{E}[\Delta_h B_t | \mathcal{P}_t^X] = \mathrm{E}[\Delta_h B_t | \mathcal{P}_t^B]$ does not converge in probability as $h \downarrow 0$, as a consequence of Proposition 10 in [4] or Proposition 13 of the current paper.

Consider expression (14). The hypothesis $E \int_0^T |b_s|\,ds < +\infty$ allows us to use the techniques of the proof of Proposition 2.5 in [6] to show that $\frac{1}{h}E[\int_t^{t+h} b_s\,ds | \mathcal{P}_t^X]$ converges in probability for almost all $t$. Using inequality (6) and the hypothesis $\mathrm{E}(|\sigma|_\alpha^p) < +\infty$, we can finally conclude that $\mathcal{P}_t^X$ is not a forward differentiating $\sigma$-field for $X$ at almost all times $t$. $\square$

4.4. *The case of fractional differential equations with analytic volatility.*

PROPOSITION 15. *Let $X \in \Xi$ be given by (8) and let $t \in (0,T)$. We assume, moreover, that $\sigma$ is a real analytic function. Then $\mathcal{P}_t^X$ is a forward differentiating $\sigma$-field for $X$ at $t$ if and only if $\sigma \equiv 0$. In this case, $\mathcal{P}^X = \{\mathcal{P}_t^X, t \in (0,T)\}$ is a discriminating collection of $\sigma$-fields and $\mathcal{P}_t^X$ is degenerate at any $t \in (0,T)$.*

PROOF. If $\sigma \equiv 0$, then $X$ is deterministic and differentiable in $t$. Consequently, $\mathcal{P}_t^X$ is a forward differentiating $\sigma$-field, but is degenerate. Assume, now, that $\sigma \not\equiv 0$. According to the Bouleau–Hirsch optimal criterion for fractional differential equations (see [10], Theorem B), we have that the law of $X_t$ is absolutely continuous with respect to the Lebesgue measure for any $t$ [indeed, we have $\text{int}\,\sigma^{-1}(\{0\}) = \varnothing$]. We deduce that $P(\sigma(X_t) = 0) = 0$ for any $t$, since $\text{Leb}(\sigma^{-1}(\{0\})) = 0$ ($\sigma$ has only isolated zeros). Proposition 14 allows us to conclude that $\mathcal{P}_t^X$ is not a forward differentiating $\sigma$-field. □

REMARK 16. The case where $\sigma$ is not assumed to be analytic seems more difficult to handle. We conjecture, however, that in this case, $\mathcal{P}_t^X$ is a forward differentiating $\sigma$-field for $X$ if and only if $t < t_x$, where $t_x$ is the deterministic time defined by

$$t_x = \inf\{t \geq 0 : x_t \notin \text{int}\,\sigma^{-1}(\{0\})\}$$

where $(x_t)_{t \in [0,T]}$ is the solution to $x_t = x_0 + \int_0^t b(x_s)\,ds$. If this conjecture is true, we would have $\xi(X) = t_x$; see (21).

## 5. Stochastic derivatives with respect to the present.

5.1. *Definition.* A consequence of Proposition 11 is that the $\sigma$-field $\mathcal{P}_t^X$ generated by $X_s$, $0 \leq s \leq t$ (the past of $X$), cannot be used for differentiating when we work with fractional Brownian motion. Moreover, we stress the following important fact: the Markov property of a Wiener diffusion $X \in \Lambda_d$ implies that taking expectations with respect to $\mathcal{P}_t^X$ produces the same effect as taking expectations only with respect to $X_t$. The following definition is then natural.

DEFINITION 17. Let $Z = (Z_t)_{t \in [0,T]}$ be a stochastic process defined on a complete probability space $(\Omega, \mathcal{F}, \mathbb{P})$ and for any $t \in (0,T)$, let $\mathcal{T}_t^Z$ be the $\sigma$-field generated by $Z_t$. We say that $Z$ admits a forward (resp., backward) stochastic derivative with respect to the present $t \in (0,T)$ if $\mathcal{T}_t^Z$ is a forward (resp., backward) differentiating $\sigma$-field for $Z$ at $t$. In this case, we set $D_+^{\mathcal{T}} Z_t := D_+^{\mathcal{T}_t^Z} Z_t$ (resp., $D_-^{\mathcal{T}} Z_t := D_-^{\mathcal{T}_t^Z} Z_t$).

EXAMPLE 18. Let $B$ be a fractional Brownian motion with Hurst index $H \in (0,1)$ and $t \in (0,T)$. Then

$$D_+^{\mathcal{T}} B_t = \begin{cases} Ht^{-1}B_t, & \text{if } H > 1/2, \\ 0, & \text{if } H = 1/2, \\ \text{does not exist}, & \text{if } H < 1/2, \end{cases}$$

and

$$D_-^{\mathcal{T}} Z_t = \begin{cases} Ht^{-1}B_t, & \text{if } H > 1/2, \\ t^{-1}B_t, & \text{if } H = 1/2, \\ \text{does not exist}, & \text{if } H < 1/2, \end{cases}$$

(see also Proposition 5). In particular, we would say that the fractional Brownian motion with Hurst index $H > 1/2$ is more regular than Brownian motion ($H = 1/2$) because of the equality between the forward and backward derivatives in the case $H > 1/2$, contrary to the case $H = 1/2$. We can identify the cause of these different regularities: the covariance function $R_H$ is differentiable along the diagonal $(t,t)$ in the case $H > 1/2$, while it is not when $H = 1/2$.

5.2. *Case of fractional differential equations.* We denote by $\Xi$ the set of fractional differential equations, that is, the subset of $\Upsilon$ whose elements are processes $X = (X_t)_{t \in [0,T]}$ solving (8) with $\sigma \in \mathcal{C}_b^2$ and $b \in \mathcal{C}_b^1$.

In the sequel, we compute $D_{\pm}^{\mathcal{T}} X_t$ for $X \in \Xi$ and $t \in (0, T)$. Let us begin with a simple case.

PROPOSITION 19. *Let $X \in \Xi$ be given by (8) and let $t \in (0, T)$. Assume, moreover, that $\sigma$ and $b$ are proportional. Then $X$ admits a forward and a backward stochastic derivative with respect to the present $t$, given by*

(22) $$D_+^{\mathcal{T}} X_t = D_-^{\mathcal{T}} X_t = Ht^{-1}\sigma(X_t)B_t + b(X_t).$$

*In particular, the present $\mathcal{T}_t^X$ is nondegenerate at $t$ if and only if $\sigma(x_0) \neq 0$ and the collection of $\sigma$-fields $\mathcal{T}^X = \{\mathcal{T}_t^X, t \in (0,T)\}$ is discriminating for $X$.*

PROOF. We will only provide the proof for $D_+^{\mathcal{T}} X_t$, the computation for $D_-^{\mathcal{T}} X_t$ being similar. Assume that $b(x) = r\sigma(x)$ with $r \in \mathbb{R}$. Then $X_t = f(B_t + rt)$, where $f : \mathbb{R} \to \mathbb{R}$ is defined by $f(0) = x_0$ and $f' = \sigma(f)$. If $\sigma(x_0) = 0$, then $X_t \equiv x_0$ and $D_+^{\mathcal{T}} X_t = 0 = \sigma(X_t)Ht^{-1}B_t + b(X_t)$. If $\sigma(x_0) \neq 0$, then it is classical that $f$ is strictly monotone. We can then write $B_t = f^{-1}(X_t) - rt$. In particular, the random variables which are measurable with respect to $X_t$ are measurable with respect to $B_t$ and vice versa. On the other hand, by using a linear Gaussian regression, it is easy to show that $D_+^{\mathcal{T}} B_t = Ht^{-1}B_t$ (see also Proposition 5). Finally, the convergences (15) and (16) and the equality (14) allow us to conclude that we have (22).

Now, let us prove that the present is nondegenerate for $X$ at $t$ if and only if $\sigma(x_0) \neq 0$. When $\sigma(x_0) = 0$, it is clear that the present is degenerate at $t$ (see the first part of this proof). On the other hand, if the present is degenerated at $t$, then there exists $c \in \mathbb{R}$ such that

$$Ht^{-1}\sigma \circ f(B_t + rt)B_t + r\sigma \circ f(B_t + rt) = c.$$

By rearranging, we obtain that $\sigma \circ f(X)(X + \alpha) = \beta$ for some $\alpha, \beta \in \mathbb{R}$, where $X = B_t + rt$. By using the fact that $X$ has a strictly positive density on $\mathbb{R}$, we deduce that $\sigma \circ f(x)(x + \alpha) = \beta$ for any $x \in \mathbb{R}$. Necessarily, $\beta = 0$ (with $x = -\alpha$) and then $f' = \sigma \circ f = 0$. We deduce that $f$ is constant and then that $f \equiv x_0$, that is, $\sigma(x_0) = 0$.

Finally, if $Ht^{-1}\sigma(X_t)B_t + b(X_t) = \sigma(X_t)(Ht^{-1}B_t + r) = 0$ a.s. for any $t$, then $\sigma(X_t) = 0 = b(X_t)$ a.s. for any $t$ and $X_t \equiv x_0$ a.s. for any $t$; see (8). In other words, the collection of $\sigma$-fields $\mathcal{T}^X = \{\mathcal{T}^X_t, t \in (0, T)\}$ is discriminating. □

Let us now describe a more general case.

THEOREM 20. *Let $X \in \Xi$ be given by (8) and let $t \in (0, T)$. Assume, moreover, that $b \in \mathcal{C}^2_b$ and that $\sigma \in \mathcal{C}^2_b$ is elliptic, that is, satisfies $\inf_{x \in \mathbb{R}} |\sigma(x)| > 0$. Then $X$ admits a forward and a backward stochastic derivative with respect to the present $t$, given by*

$$
\begin{aligned}
D^{\mathcal{T}}_+ X_t &= D^{\mathcal{T}}_- X_t \\
&= b(X_t) + H\frac{\sigma(X_t)}{t}\bigg\{\int_0^{X_t} \frac{dy}{\sigma(y)} \\
&\qquad - \mathrm{E}\bigg[\int_0^t \frac{b}{\sigma}(X_s)\,ds + \int_0^t \int_0^t \beta^H_r(s)\,\delta B_r\,ds \\
&\qquad\qquad - t\int_0^t \beta^H_r(t)\delta B_r \Big| X_t\bigg]\bigg\},
\end{aligned}
$$
(23)

*where*

$$\beta^H_r(t) = \bigg(\mathcal{O}_H \int_0^r \frac{b'\sigma - b\sigma'}{\sigma}(X_s)\mathbf{1}_{s\geq .}\,ds\bigg)(t).$$

*Recall that $\mathcal{O}_H$ is defined by (4).*

PROOF. First, note that $\beta^H_r(t)$ belongs to the domain of the divergence operator $\delta^B$, due to the additional hypothesis on $b$ and $\sigma$. We provide only the proof for $D^{\mathcal{T}}_+ X_t$, the computation for $D^{\mathcal{T}}_- X_t$ being similar. Földmer [6], Section 4, tackles the problem of the computation of the time reversed drift of a non-Markovian diffusion by means of a Girsanov transformation and the Malliavin calculus. Our proof uses such a strategy, coupled with the transfer principle.

*First step.* Assume that $\sigma \equiv 1$. Using the transfer principle and the isometry $\mathcal{K}_H$, it holds that

$$X_t = \int_0^t K_H(t,s) \, dY_s,$$

where

$$Y_t = W_t + \int_0^t a_r \, dr.$$

Here, we set

$$a_r = \left(\mathcal{K}_H^{-1} \int_0^\cdot b(X_s) \, ds\right)(r).$$

We know (see [13], Theorem 2) that the process $X = (X_t)_{t \in [0,T]}$ is a fractional Brownian motion under the new probability measure $\mathbb{Q} = G \cdot \mathbb{P}$, where

$$G = \exp\left(-\int_0^T a_s \, dW_s - \tfrac{1}{2} \int_0^T a_s^2 \, ds\right).$$

Using the integration by parts, of Malliavin calculus, we can write, for $g : \mathbb{R} \to \mathbb{R} \in C_b^1$,

$$\begin{aligned}
\mathrm{E}[(X_{t+h} - X_t)g(X_t)] &= \mathrm{E}_{\mathbb{Q}}[G^{-1} g(X_t) \delta^X(\mathbf{1}_{[t,t+h]})] \\
&= \mathrm{E}_{\mathbb{Q}}[G^{-1} \langle \mathbf{1}_{[t,t+h]}, D^X g(X_t) \rangle_{\mathcal{H}}] \\
&\quad + \mathrm{E}_{\mathbb{Q}}[g(X_t) \langle \mathbf{1}_{[t,t+h]}, D^X G^{-1} \rangle_{\mathcal{H}}] \\
&= \mathrm{E}[g'(X_t)] \langle \mathbf{1}_{[t,t+h]}, \mathbf{1}_{[0,t]} \rangle_{\mathcal{H}} \\
&\quad + \mathrm{E}[G g(X_t) \langle \mathcal{K}_H^* \mathbf{1}_{[t,t+h]}, \mathcal{K}_H^* D^X G^{-1} \rangle_{L^2}].
\end{aligned}$$

But $\mathcal{K}_H^* D^X G^{-1} = D^Y G^{-1}$ (transfer principle). Since

$$G^{-1} = \exp\left(\int_0^T a_s \, dY_s - \tfrac{1}{2} \int_0^T a_s^2 \, ds\right),$$

we have

$$\begin{aligned}
G \times D_t^Y(G^{-1}) &= a_t + \int_0^T D_t^Y a_s \, dY_s - \int_0^T a_s D_t^Y a_s \, ds \\
&= a_t + \int_0^T D_t^Y a_s \, dW_s.
\end{aligned}$$

Moreover,

$$\int_0^T D_s^Y a_r \, dW_r = \int_0^T (\mathcal{K}_H^* D_s^X a)(r) \, dW_r = \int_0^T D_s^X a_r \, \delta B_r := \Phi(s)$$

and
$$(\mathcal{K}_H^* \mathbf{1}_{[0,t]})(s) = K_H(t,s)\mathbf{1}_{[0,t]}(s).$$

Therefore,
$$\langle \mathcal{K}_H^* \mathbf{1}_{[t,t+h]}, G\mathcal{K}_H^* D^X G^{-1} \rangle_{L^2}$$
$$= (\mathcal{K}_H a)(t+h) - (\mathcal{K}_H a)(t) + (\mathcal{K}_H \Phi)(t+h) - (\mathcal{K}_H \Phi)(t)$$
$$= \int_t^{t+h} b(X_u)\, du + (\mathcal{K}_H \Phi)(t+h) - (\mathcal{K}_H \Phi)(t).$$

By the stochastic Fubini theorem, we have $(\mathcal{O}_H \Phi)(t) = \int_0^T (\mathcal{O}_H D_\cdot^X a_r)(t)\, \delta B_r$. We set
$$\beta_r^H(t) = (\mathcal{O}_H D_\cdot^X a_r)(t) = \left( \mathcal{O}_H \int_0^r b'(X_s) \mathbf{1}_{s\geq \cdot}\, ds \right)(t).$$

We then deduce that
$$\mathrm{E}[(X_{t+h} - X_t)g(X_t)]$$
(24)
$$= \mathrm{E}[g'(X_t)]\langle \mathbf{1}_{[t,t+h]}, \mathbf{1}_{[0,t]} \rangle_{\mathcal{H}}$$
$$+ \mathrm{E}\left[ g(X_t)\left( \int_t^{t+h} b(X_s)\, ds + \int_t^{t+h} \int_0^T \beta_r^H(s)\, \delta B_r\, ds \right) \right].$$

By developing $\mathrm{E}[X_t g(X_t)]$ as in (24), we obtain
$$t^{2H} \mathrm{E}[g'(X_t)] = \mathrm{E}\left[ g(X_t)\left( X_t - \int_0^t b(X_s)\, ds - \int_0^t \int_0^T \beta_r^H(s)\, \delta B_r\, ds \right) \right].$$

Then
$$\mathrm{E}[\Delta_h X_t \mid X_t]$$
$$= h^{-1} \langle \mathbf{1}_{[t,t+h]}, \mathbf{1}_{[0,t]} \rangle_{\mathcal{H}} \left( X_t - \mathrm{E}\left[ \int_0^t b(X_s)\, ds + \int_0^t \int_0^T \beta_r^H(s)\, \delta B_r\, ds \Big| X_t \right] \right)$$
$$+ h^{-1} \mathrm{E}\left[ \int_t^{t+h} b(X_s)\, ds + \int_t^{t+h} \int_0^T \beta_r^H(s)\, \delta B_r\, ds \Big| X_t \right].$$

We deduce that $\mathrm{E}[\Delta_h X_t \mid X_t]$ converges in probability, as $h \downarrow 0$, to
$$b(X_t) + \frac{H}{t} X_t - \frac{H}{t} \mathrm{E}\left[ \int_0^t b(X_s)\, ds + \int_0^t \int_0^T \beta_r^H(s)\, \delta B_r\, ds - \int_0^T \beta_r^H(t)\, \delta B_r \Big| X_t \right].$$

Since $\lim_{h \downarrow 0} \mathrm{E}[\Delta_h X_t \mid X_t]$ does not depend on $T$, we finally obtain (23) in the particular case where $\sigma \equiv 1$ by letting $T \downarrow t$.

*Second step.* Assume that $\sigma$ does not vanish. Set $Y_t = h(X_t)$, where $h(x) = \int_0^x \frac{dy}{\sigma(y)}$. Using the change of variables formula, we obtain that $Y$ satisfies
$$Y_t = y_0 + B_t + \int_0^t \frac{b}{\sigma} \circ h^{-1}(Y_s)\, ds, \qquad t \in [0,T].$$

Since, on the one hand, the $\sigma$-fields generated by $X_t$ and $Y_t$ are the same and, on the other hand, $X$ has $\alpha$-Hölder continuous paths with $\alpha > 1/2$, we have

$$D_+^{\mathcal{T}} X_t = \sigma(X_t) D_+^{\mathcal{T}} Y_t.$$

Expression (23) is then a consequence of the first step of the proof. $\square$

REMARK 21. When $\sigma$ does not vanish and $b \equiv r\sigma$ with $r \in \mathbb{R}$, we can apply either Proposition 19 or Theorem 20 to compute $D_\pm^{\mathcal{T}} X_t$. Of course, the conclusions are the same. Indeed, since we have, in this case, $b'\sigma - b\sigma' \equiv 0$ and $\int_0^{X_t} \frac{dy}{\sigma(y)} = B_t + rt$ [since $X_t = f(B_t + rt)$ with $f$ satisfying $f' = \sigma \circ f$], formula (23) can be simplified to (22).

Compared to the case where $\sigma$ and $b$ are proportional, here, it is more difficult to decide if the present (i.e., the collection of $\sigma$-fields generated by $X_t$) is discriminating or not.

In the framework of the stochastic embedding of dynamical systems introduced in [2], the set of processes, called *Nelson differentiable processes*, which satisfy the equality between stochastic forward and stochastic backward derivatives plays a fundamental role (see [3], Chapters 3 and 7). We stress on the fact that solutions of stochastic differential equations driven by a fractional Brownian motion with $H > 1/2$ provide examples of Nelson differentiable processes which are not absolutely continuous.

**Acknowledgments.** This work was carried out while the first author benefited from the hospitality of Université Paris 6. This institution is gratefully acknowledged. We also wish to thank the anonymous referee for a careful and thorough reading of this work and for his constructive remarks. We finally acknowledge C. Stricker for pointing out an error in a previous version.

LABORATOIRE DE MATHÉMATIQUES
UNIVERSITÉ DE FRANCHE-COMTÉ
16 ROUTE DE GRAY
25030 BESANCON CEDEX
FRANCE
E-MAIL: sedarses@ccr.jussieu.fr

LPMA
UNIVERSITÉ PARIS 6
BOÎTE COURRIER 188
75252 PARIS CEDEX 05
FRANCE
E-MAIL: nourdin@ccr.jussieu.fr